\documentclass[conference]{IEEEtran}
\IEEEoverridecommandlockouts
\usepackage{cite}
\usepackage{amsmath,amssymb,amsfonts}
\usepackage{algorithmic}
\usepackage{graphicx}
\usepackage{textcomp}
\usepackage{xcolor}

\usepackage{physics}
\usepackage{caption}
\usepackage{url}
\usepackage{balance}
\usepackage{subfig}
\usepackage{hyperref}

\renewcommand{\vec}{\boldsymbol}

\graphicspath{{./figures/}}

\def\BibTeX{{\rm B\kern-.05em{\sc i\kern-.025em b}\kern-.08em
		T\kern-.1667em\lower.7ex\hbox{E}\kern-.125emX}}
\begin{document}
	
	\title{Dynamic mode decomposition as an analysis tool for time-dependent 
	partial differential equations}
	
	\author{
		\IEEEauthorblockN{Miha 
		Rot\IEEEauthorrefmark{1}\IEEEauthorrefmark{2}, Martin 
		Horvat\IEEEauthorrefmark{3}, Gregor Kosec\IEEEauthorrefmark{2}}
		
		\IEEEauthorblockA{\IEEEauthorrefmark{1}Jozef Stefan International 
			Postgraduate School, Jamova cesta 39, 1000 Ljubljana, Slovenia}
		\IEEEauthorblockA{\IEEEauthorrefmark{2}``Jožef Stefan'' Institute, 
			Parallel and Distributed Systems Laboratory, Jamova cesta 39, 1000 
			Ljubljana, Slovenia}
		\IEEEauthorblockA{\IEEEauthorrefmark{3}University of Ljubljana, Faculty 
			of Mathematics and Physics, Department of Physics,  Jadranska cesta 
			19, SI-1000 Ljubljana. \\
			miha.rot@ijs.si}
		
		\thanks{The authors would like to acknowledge the financial
			support of the Slovenian Research Agency (ARRS) research core 
			funding 
			No.\ P2-0095 and the Young Researcher program PR-10468.}
	}
	
	\maketitle

\begin{abstract}
	The time-dependent fields obtained by solving partial differential 
	equations in two and more dimensions quickly overwhelm the analytical 
	capabilities of the human brain. A meaningful insight into the temporal 
	behaviour can be obtained by using scalar reductions, which, however, come 
	with a loss of spatial detail. Dynamic Mode Decomposition is a data-driven 
	analysis method that solves this problem by identifying oscillating spatial 
	structures and their corresponding frequencies. This paper presents the 
	algorithm and provides a physical interpretation of the results by applying 
	the decomposition method to a series of increasingly complex examples.
\end{abstract}

\bigskip

\begin{IEEEkeywords}
	\textit{\textbf{dynamic mode decomposition; data-driven analysis; wave 
	equation; non-Newtonian natural convection}}
\end{IEEEkeywords}


\section{Introduction}

Humans posses a remarkable ability for pattern recognition that is however 
still quickly overwhelmed when analysing complex systems. When dealing with 
time series of two or three dimensional field snapshots, as is common in e.g. 
hydrodynamics, we are often forced to use reductions that decrease complexity 
and provide an approachable insight into the system's behaviour. An example of 
such reduction when dealing with natural convection is tracking the ratio 
between heat convection and conduction, also known as the Nusselt number, on a 
significant boundary. The Nusselt number provides us with a scalar value that 
reflects the behaviour in our system. Analysing time series of such values can 
provide insight on whether the system is stable, oscillatory or maybe even 
chaotic, with Fourier analysis available as a potent tool to extract the 
characteristic frequencies reflected in the chosen reduction. No matter how 
good the reduction is we still lose the majority of the spatial detail which is 
often important for understanding the dynamics.

Spatial structure of the analysed system is required for good understanding of 
the dynamics, which is especially important when dealing with engineering 
challenges like the airflow around aircraft structures \cite{Ricciardi2019}. A 
common approach to identifying important spatial structures in analysed fields 
is modal decomposition, with Proper Orthogonal Decomposition 
(POD)\footnote{Also known as principal component analysis (PCA) and various 	
other names in different fields.}, proposed for hydrodynamics by Lumley 
\cite{Lumley1967} in 1967, and it's variations \cite{Rowley2005, Taira2017} as 
one of the most widely used techniques. POD is a decomposition technique that 
identifies orthogonal modes that best\footnote{Optimally in $L_2$ sense that 	
corresponds to energy when decomposing flow field \cite{Holmes2012}.} represent 
the dataset and provides an ordering based on their importance.

Modes obtained with POD provide information about energetically important parts 
of the system but are completely oblivious to the temporal behaviour meaning 
that the results of the decomposition would stay the same even if analysis was 
performed on reordered snapshots. Dynamic mode decomposition (DMD) was proposed 
by Schmid in a 2008 talk and the subsequent paper in 2010 \cite{Schmid2010} as 
a technique that joins the spatial aspects of POD and temporal aspects of 
Fourier transform \cite{Taira2017}. DMD identifies characteristic frequencies, 
corresponding spatial structures and whether they amplify, decay or remain 
constant through the sampled timespan. The method was initially proposed for 
hydrodynamics but has since found other applications including infectious 
disease spread \cite{Proctor2015} and computer vision \cite{Grosek2014} as it 
is completely data-driven and independent of the underlying dynamics.

This paper introduces the DMD algorithm and interpretation in Sec. 
\ref{ch:DMD}, demonstrates the results on vibrating membrane examples in Sec. 
\ref{ch:membrane} and concludes with the decomposition applied to a more 
interesting case of oscillatory non-Newtonian natural convection in Sec. 
\ref{ch:DVD}.

\section{Dynamic Mode Decomposition}
\label{ch:DMD}

\subsection{Interpretation}
The DMD algorithm is based on finding eigenvalues $\lambda$ and eigenvectors 
$\vec{\varphi}$ of the linear mapping represented by the matrix $\vec{A}$
\begin{equation}
	\vec{v}_{i+1} = \vec{A} \vec{v}_i,
	\label{eq:linearMap}
\end{equation}
that connects the subsequent states $\vec{v}_i$ and $\vec{v}_{i+1}$ of the 
analysed system. As such the algorithm is completely independent of the 
governing equations that drive the underlying dynamics and can be applied to 
any experimental or simulated data.

The linear nature of eigendecomposition does not preclude us from analysing 
non-linear systems as the dynamics of any such system can be expressed with an 
infinite-dimensional linear Koopman operator \cite{Rowley2009} with the DMD 
providing a good approximation of it's eigendecomposition as long as the 
quality and quantity of input data is sufficient \cite{Tu2014}.

DMD eigenvalues $\lambda$ provide information about the temporal behaviour of 
the spatial structure described in the eigenvector. Eigenvalues are complex 
numbers with the mode's angular frequency information contained in the argument
\begin{equation}
	\label{eq:frequency}
	\nu_i = \frac{\arg(\lambda_i)}{2\pi \Delta t},
\end{equation}
with sampling interval $\Delta t$ providing scaling to the output, that is 
otherwise unaffected by the frequency of data snapshot sampling. The magnitude 
of eigenvalues provides information about changes in mode's strength, 
with $|\lambda| > 1$ for amplifying and $|\lambda| < 1$ for decaying modes. 

The linear combination of DMD modes calculated from real valued input data has 
to be real valued, which can be satisfied either by real 
eigenvalues/eigenvectors or complex-conjugate pairs of complex 
eigenvalues/eigenvectors. Real eigenvalues represent modes with frequency 0. At 
least one such background mode with frequency 0 and eigenvalue magnitude of 1 
is present whenever the decomposed data has non-zero average but multiple such 
real-valued growing/decaying non-oscillatory modes can appear. Complex 
conjugate pairs represent oscillating modes and can be treated as a single mode 
and will be in most of the later visualisations.

\subsection{Algorithm}
The algorithm for exact DMD proposed by Tu \cite{Tu2013thesis} starts by 
taking subsequent snapshots $\vec{v}_i$ of analysed fields separated by 
$\Delta t$. Uniformly sampled snapshots simplify calculation and frequency 
interpretation but are not explicitly required for DMD. Snapshots are 
flattened from their arbitrary shape into vectors with length $M$ and arranged 
as columns in a matrix
\begin{equation}
	\label{eq:vector}
	\vec{V}_i^N = \mqty[\vec{v}_i & \vec{v}_{i+1} & \cdots & \vec{v}_N],
\end{equation}
with $N$ denoting the number of snapshots. The number of data points in a 
snapshot $M$ is usually much larger than the number of snapshots $N$, leading 
to tall and narrow matrices $\vec{V}_i^N$ that can be exploited to efficiently 
calculate the eigen-decomposition of the large $M \times M$ matrix $\vec{A}$.

The snapshot matrix $\vec{V}_i^N$ needs to have a sufficient\footnote{Generally 
	at least two times the number of relevant oscillatory modes.} rank to 
	describe 
the dynamics. Rank deficiency issues usually appear when dealing with standing 
waves \cite{Tu2014, Kutz2016} that cause a linear dependency between snapshots 
in columns of $\vec{V}_i^N$ and can be solved by stacking subsequent snapshots 
\begin{equation}
	\vec{V}_i^N = \mqty[\vec{v}_i & \vec{v}_{i+1} & \cdots & \vec{v}_{N-m} \\ 
	\vec{v}_{i+1} & \vec{v}_{i+2} & \cdots & \vec{v}_{N-m+1} \\ \vdots & \vdots 
	& \ddots & \vdots \\
	\vec{v}_{i + m} & \vec{v}_{i + m + 1} & \cdots & \vec{v}_{N}],
	\label{eq:stacking}
\end{equation}
ensuring $\rank(\vec{V}_i^N)>=m$ at the cost of effectively larger system, and 
only using the first $M$ values from the now larger DMD eigenvectors. Time 
shift augmentation can also help with noisy data or in cases where snapshots 
are small \cite{Kutz2016}.

For simplicity we continue our discussion with the notation from 
Eq.~\eqref{eq:vector}, apply the relationship from Eq. \eqref{eq:linearMap} 
to the matrix
\begin{equation}
	\vec{V}_1^N = \mqty[\vec{A}\vec{v}_0 & \vec{A}\vec{v}_{i} & \cdots &
	\vec{A}\vec{v}_{N-1}] = \vec{A}\vec{V}_0^{N-1},
\end{equation}
and perform the compact singular value decomposition (SVD) on the data matrix 
\begin{equation}
	\vec{V}_0^{N-1} = \vec{U} \vec{S} \vec{W}^\dag.
\end{equation}
The results of the SVD can already be assigned physical meaning as it is one of 
the two main algorithms used to compute POD. Singular values denoting the 
importance of POD modes are contained in the diagonal $\vec{S}$ with the
corresponding eigenvectors in columns of $\vec{U}$. The complexity of the 
underlying dynamics can be estimated from the spectrum of singular values. 
Cases with exponentially decreasing singular value magnitudes signify a 
low-dimensional behaviour that can be described with a handful of modes. Slow 
decay is not as simple to explain as it can be a product of either very 
complex dynamics or noisy data.

Output of the SVD algorithm can be truncated based on singular values or prior 
knowledge about the system to retain only $r$ singular values and corresponding 
singular vectors in $\vec{U}$ and $\vec{W}$. Truncation is beneficial as it 
reduces the number of numerical operations, as evident from the following 
paragraph, and removes noise, but care has to be taken in order to avoid 
truncating relevant parts of the dynamics. 
More details about truncation with examples and advanced strategies are 
available in \cite{Kutz2016}.

The optionally truncated result can be reorganized and rewritten as
\begin{equation}
	\vec{U}^\dag \vec{V}_1^N \vec{W} \vec{S}^{-1} = \vec{U}^\dag \vec{A} 
	\vec{U} := \widetilde{\vec{A}},
\end{equation}
with $\widetilde{\vec{A}}$ representing the projection of $\vec{A}$ onto 
the POD eigenvectors. We then perform the eigendecomposition of 
$\widetilde{\vec{A}}$ which is a much smaller $r \times r$ matrix than the 
full, $M \times M$ sized, matrix $\vec{A}$. The decomposition
\begin{equation}
	\widetilde{\vec{A}} \vec{w} = \lambda \vec{w}.
\end{equation}
gives us eigenvalues $\lambda_i$ and eigenvectors $\vec{w}_i$. The eigenvalues 
are the same for both $\widetilde{\vec{A}}$ and $\vec{A}$ while the 
eigenvectors $\vec{\varphi}$ for the latter still need to be transformed 
from their projected representation
\begin{equation}
	\vec{\varphi} = \frac{1}{\lambda} \vec{V}_1^N  \vec{W} \vec{S}^{-1} 
	\vec{w}.
\end{equation}

\subsection{Mode ordering}
\label{sec:ordering}
The main disadvantage of DMD is that the hierarchy of the identified modes 
is unclear. Their relative importance can not be deduced from either the 
unitary eigenvectors or the eigenvalues that carry temporal data. This 
deficiency is bypassed by calculating the projection coefficients for 
expressing the measured datasets in terms of DMD eigenvectors
\begin{equation}
	\vec{v}(t) = \sum_k \vec{\varphi}_k b_k(t).
\end{equation}
As the relative importance of eigenvectors in an oscillating dataset 
changes with time it is best to calculate their importance for all 
available snapshots expressed in matrix form as
\begin{equation}
	\vec{V}_0^{N-1} = \vec{\Phi} \vec{B},
\end{equation}
with eigenvectors collected in columns of matrix $\vec{\Phi} = \mqty[ 
\vec{\varphi}_0 & \vec{\varphi}_1 & \cdots & \vec{\varphi}_{r} ]$ and 
coefficients corresponding to datasets in $\vec{V}_0^{N-1}$ 
collected as columns in matrix $\vec{B} = \mqty[ \vec{b}_0 & \vec{b}_1 & 
\cdots, \vec{b}_{N-1} ]$. The complex number $\vec{B}_{ij}$ corresponds to 
the projection coefficient of the $i$-th eigenvector for the $j$-th data 
sample. Columns of $\vec{B}$ can be thus interpreted as the time evolution 
of mode importance that can be condensed into a single value by 
calculating the average value. We use the mode power ($P$)
\begin{equation}
	P_i = \frac{1}{N} \sum_{j=0}^{N-1} \abs{B_{ij}}
	\label{eq:modePower}
\end{equation}
to rank the calculated modes.

An efficient algorithm for calculating $\vec{B}$ has also been proposed by Tu 
et al. \cite{Tu2014}. The algorithm avoids additional computation beyond 
calculating the left eigenvectors $\vec{Z}$ during the eigendecomposition 
of $\widetilde{\vec{A}}$ by expressing the matrix as
\begin{equation}
	\vec{B} = \vec{Z}^\dag \vec{S} \vec{W}^\dag.
\end{equation}

\section{Vibrating membrane}
\label{ch:membrane}

\subsection{Numerical model}
We demonstrate the DMD algorithm with a damped wave equation on a 2D square 
membrane. This case is convenient as it still allows for comparison against 
known analytic results but is closer to possible real world applications. We 
are solving the wave equation 
\begin{equation}
	\label{eq:waveEquation}
	\pdv[2]{h}{t} = c^2 \laplacian h - \gamma \pdv{h}{t},
\end{equation}
with $c$ denoting the velocity of wave propagation and $\gamma$ the damping 
factor, on a square domain $\Omega = \{(x, y) : x, y \in [0, 1] \}$ with 
Dirichlet boundary condition $h_{\partial\Omega} = 0$. A constant velocity 
$c=1$ is used for all of the presented examples.

We start the dynamics with a Gaussian initial condition with $\sigma = 0.1$ and 
amplitude of $1$ placed at an arbitrarily chosen position $(0.3, 0.4)$ to 
simulate an impact against the membrane. We would like to show that DMD 
identifies eigenmodes for the membrane and analyse how damping might affect the 
accuracy.

\subsection{Analytic solution}

\begin{figure}
	\includegraphics[width=\linewidth]{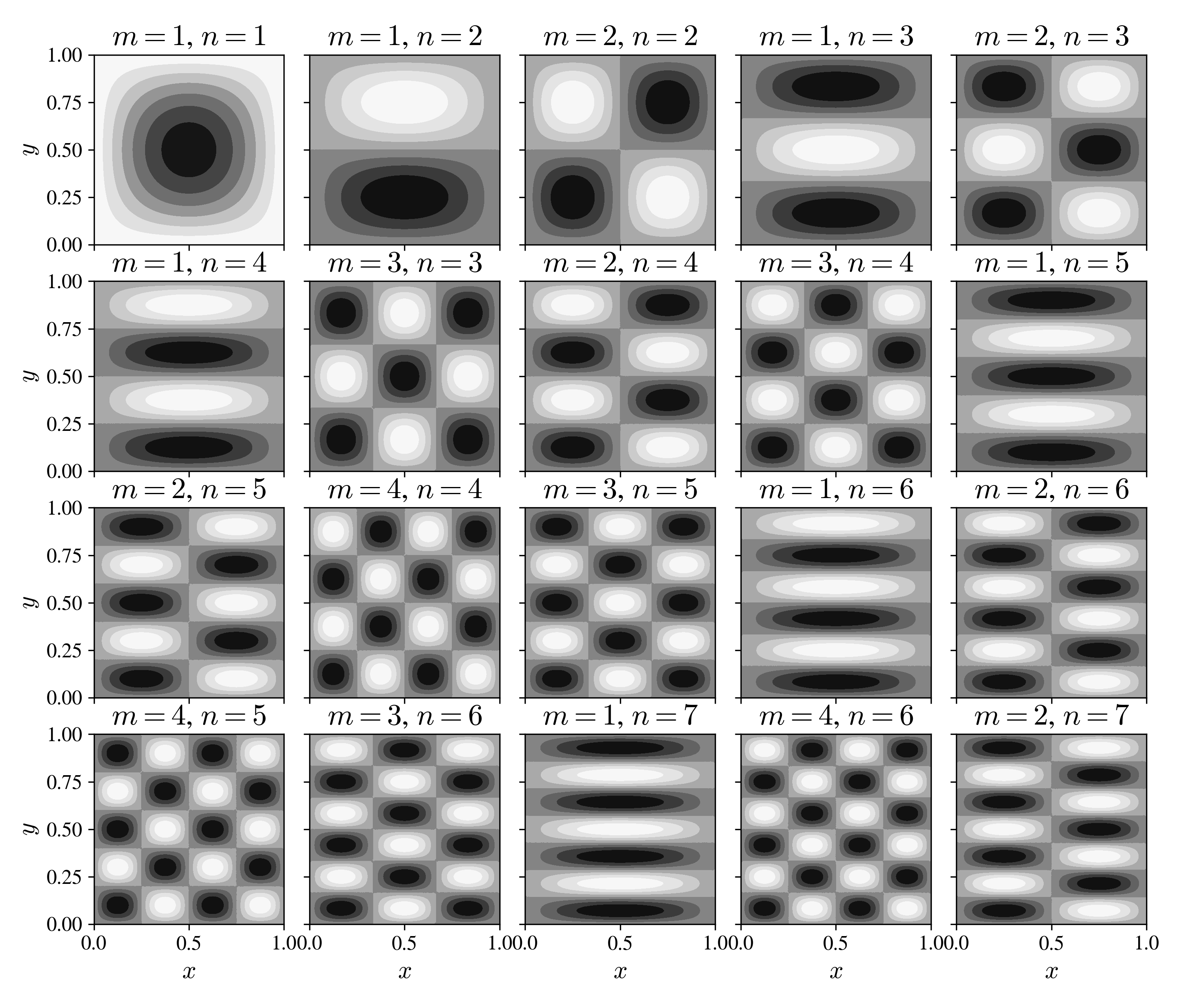}
	\caption{Analytic eigenvectors for the oscillation of a square 
		membrane, ordered by increasing frequency.}
	\label{fig:rectangularVectorsAnalytic}
\end{figure}

We can obtain an analytic solution for the partial differential Eq. 
\eqref{eq:waveEquation} by using the separation of variables. For now we ignore 
damping as it is not relevant for the eigenstate spatial configuration. 
The solution is a linear combination of modes defined with two integer 
quantization indices $m, n >=1$. Analytical eigenfrequencies and eigenvectors 
are expressed in terms of $m$ and $n$
\begin{equation}
	\nu = \frac{c}{2} \sqrt{m^2 + n^2},
\end{equation}
\begin{equation}
	h(x, y, m, n) = \sin(m \pi x) \sin(n \pi y),
\end{equation}
with eigenvectors shown in Figure~\ref{fig:rectangularVectorsAnalytic}. Only 
one eigenvector is shown for the degenerate frequencies.

\subsection{Results}

\begin{figure}
	\includegraphics[width=\linewidth]{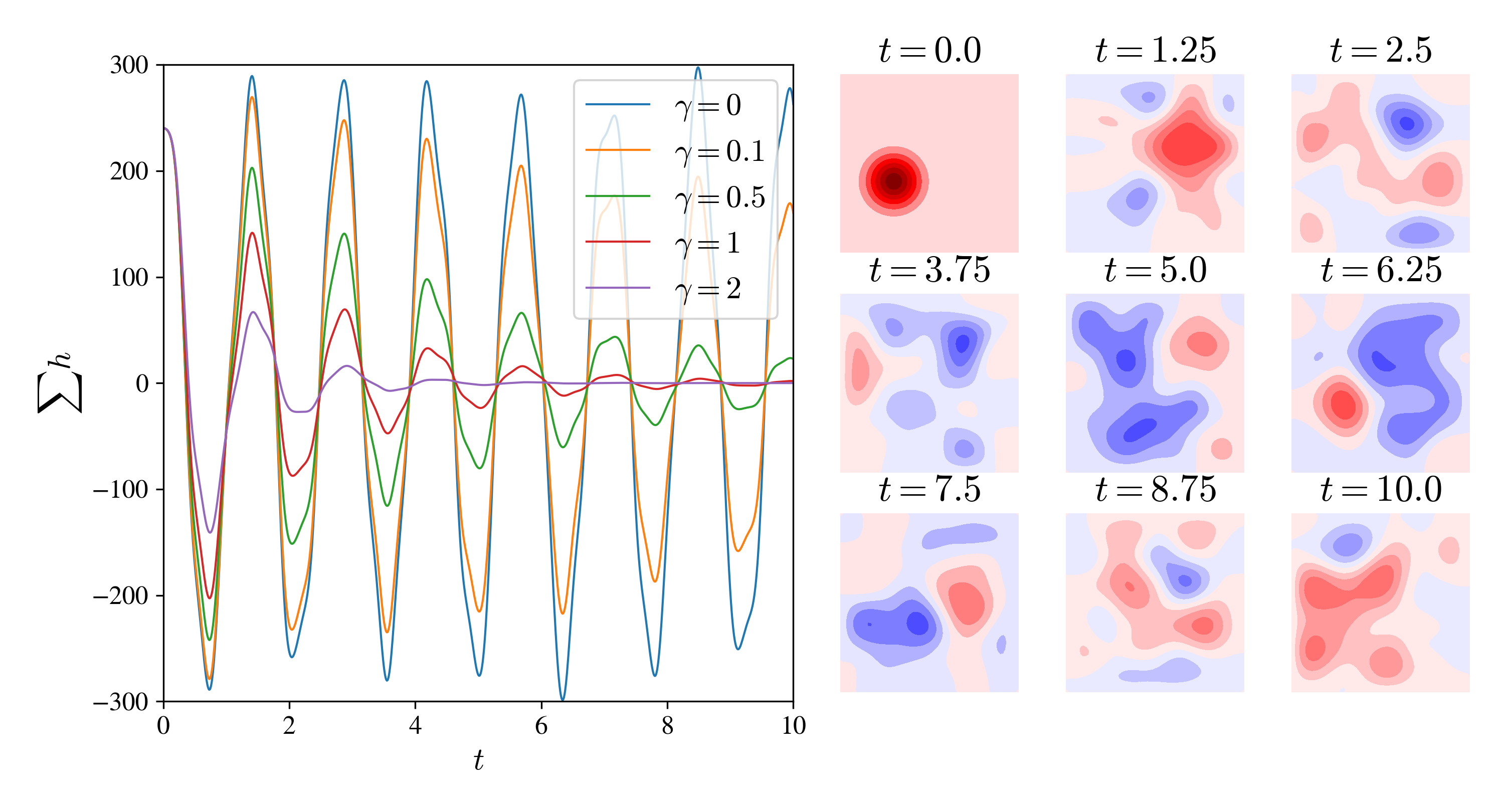}
	\caption{\textit{left:} The sum of height values provides a  
		simple insight into the system's behaviour. The time evolution of thus 
		derived observable is shown for different damping rates $\gamma$.
		\textit{right:} A selection of snapshots used in the decomposition 
		providing an insight into the system's behaviour during the observed 
		time. A heat map display is used to visualise the scalar values with 
		blues for negative, white for zero and reds for positive values.}
	\label{fig:rectangularSnapshots}
\end{figure}

The partial differential equation is solved numerically for different damping 
rates using the meshless RBF-FD~\cite{tolstykh2003rbffd} method implemented 
with the Medusa \cite{medusa} C++ library. A subset of system state snapshots 
for the 0 damping case that that we use in DMD is shown in the right graphs of 
Figure~\ref{fig:rectangularSnapshots}. The dynamics are disordered due to the 
non-symmetric position of the initial jolt and do not offer much insight at the 
first glance. We can use the sum of membrane heights as an observable that 
reduces the system's behaviour into a single scalar value to visualise the 
oscillation and the effect of different damping rates as shown in the left 
graph of Figure~\ref{fig:rectangularSnapshots}. 

We use DMD on 1000 uniformly sampled snapshots from the interval $t \in [0, 
10]$ that are composed into the data matrix $\vec{V}$ with stacking $m=10$. 
Each snapshot is composed of 3985 height values in computational nodes for the 
mesh-free partial differential equation solving, which are uniformly 
positioned within $\Omega$.

\begin{figure}
	\includegraphics[width=\linewidth]{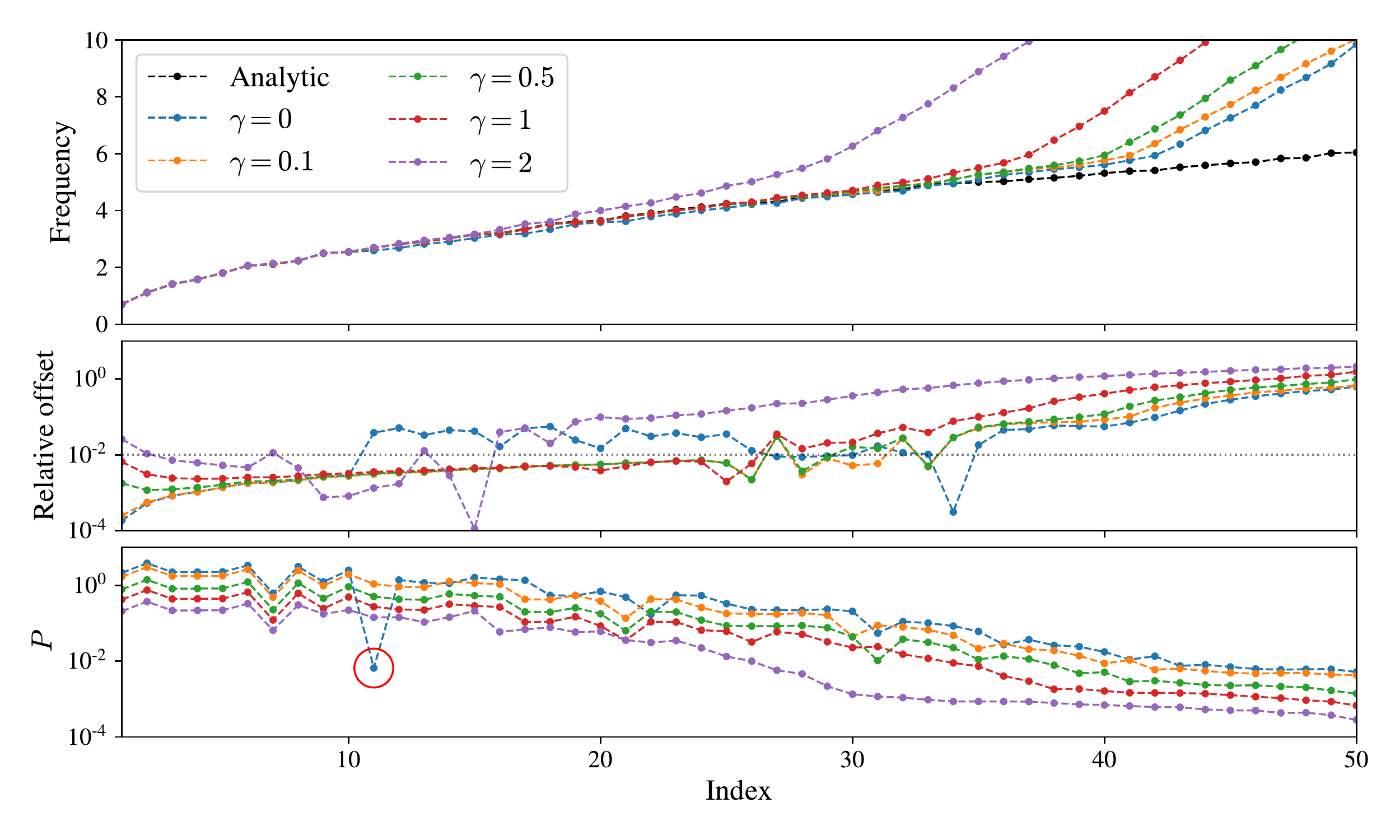}
	\caption{\textit{top:} Eigenfrequencies calculated with DMD for different 
	damping rates compared to analytic values. \textit{centre:} The relative 
	error between DMD eigenfrequencies and analytic values for the square 
	membrane. \textit{bottom:} Mode powers for the calculated DMD modes.}
	\label{fig:rectangularFrequencyComparison}
\end{figure}

The resulting frequencies, calculated from DMD eigenvalues as described in 
Eq.~\eqref{eq:frequency}, are shown in the top graph of
Figure~\ref{fig:rectangularFrequencyComparison}. The results provide a good 
estimation for the analytic results with sub percent relative error for the 
first $\sim 30$ modes as shown in the central graph. Estimated frequencies 
remain good even for the relatively strongly damped $\gamma=1$ case showing 
promise for eigenfrequency estimation on experimental data obtained from 
complex objects.

There is a slight inconsistency in the results with a sharp increase in the 
relative error for the $\gamma = 0$ frequencies. We would expect that the 
least damped case would give the best results but there seems to be an 
additional non-physical mode at index 11 that has caused a shift in the 
following frequencies. Fortunately we can use the mode power described in 
Sec.~\ref{sec:ordering} and shown in the bottom graph of 
Figure~\ref{fig:rectangularFrequencyComparison} to detect and remove the 
highlighted spurious mode with a significantly lower strength. It is surprising 
how little nonsensical modes we identified even though we initially disturbed 
the membrane in a non-symmetric location and how chaotic the actual system 
shown in the heat maps of Figure~\ref{fig:rectangularSnapshots} looks to a 
human observer. 

\begin{figure}
	\includegraphics[width=\linewidth]{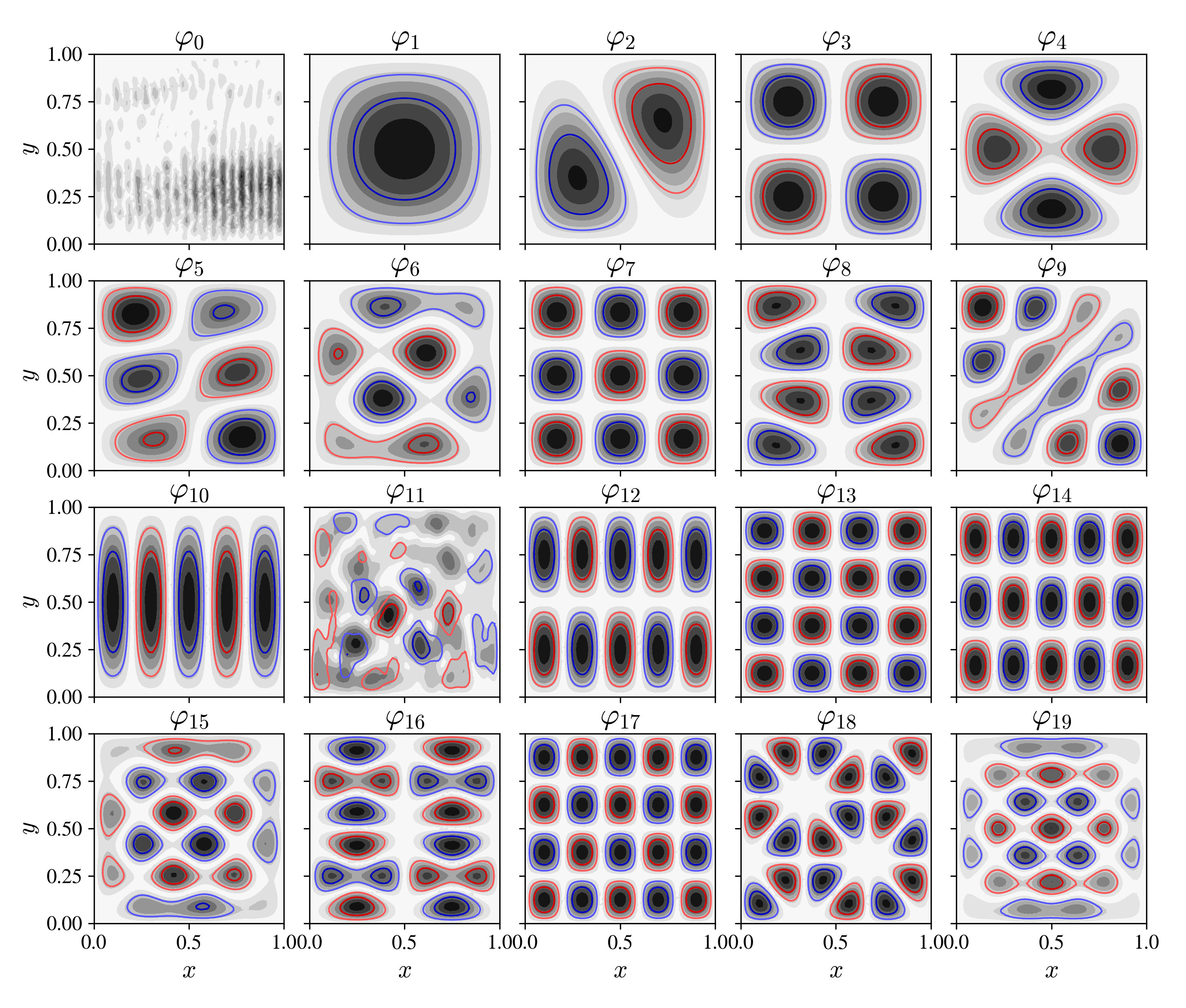}
	\caption{DMD mode eigenvectors ordered by frequency for the $\gamma = 0$ 
	case. The positions of modes are shifted by one relative to 		
	Figure~\ref{fig:rectangularFrequencyComparison} because of the constant, 
	zero frequency, background mode. The grayscale image displays absolute 	
	value with contours ranging from blue to red for the imaginary part.}
	\label{fig:rectangularVectorsDMD}
\end{figure}

We can also compare the analytic eigenvectors shown in 
Figure~\ref{fig:rectangularVectorsAnalytic} and DMD eigenvectors in 
Figure~\ref{fig:rectangularVectorsDMD}. The DMD vector positions in the grid 
are shifted by one compared to the analytic due to the non-oscillatory 
background mode that is always present in the decomposition of data with 
non-zero average. Some eigenvectors match almost perfectly while others wildly 
differ. This might seem inconsistent at the first glance but can be explained 
with degeneracy in the analytical spectrum. Linear combinations of degenerate 
eigenvectors still constitute a valid eigenvector. The non-symmetric 
eigenvector $\varphi_{11}$ can quickly be identified as the one corresponding 
to the previously mentioned spurious mode.

\subsection{Oddly shaped membrane}
\begin{figure}
	\includegraphics[width=\linewidth]{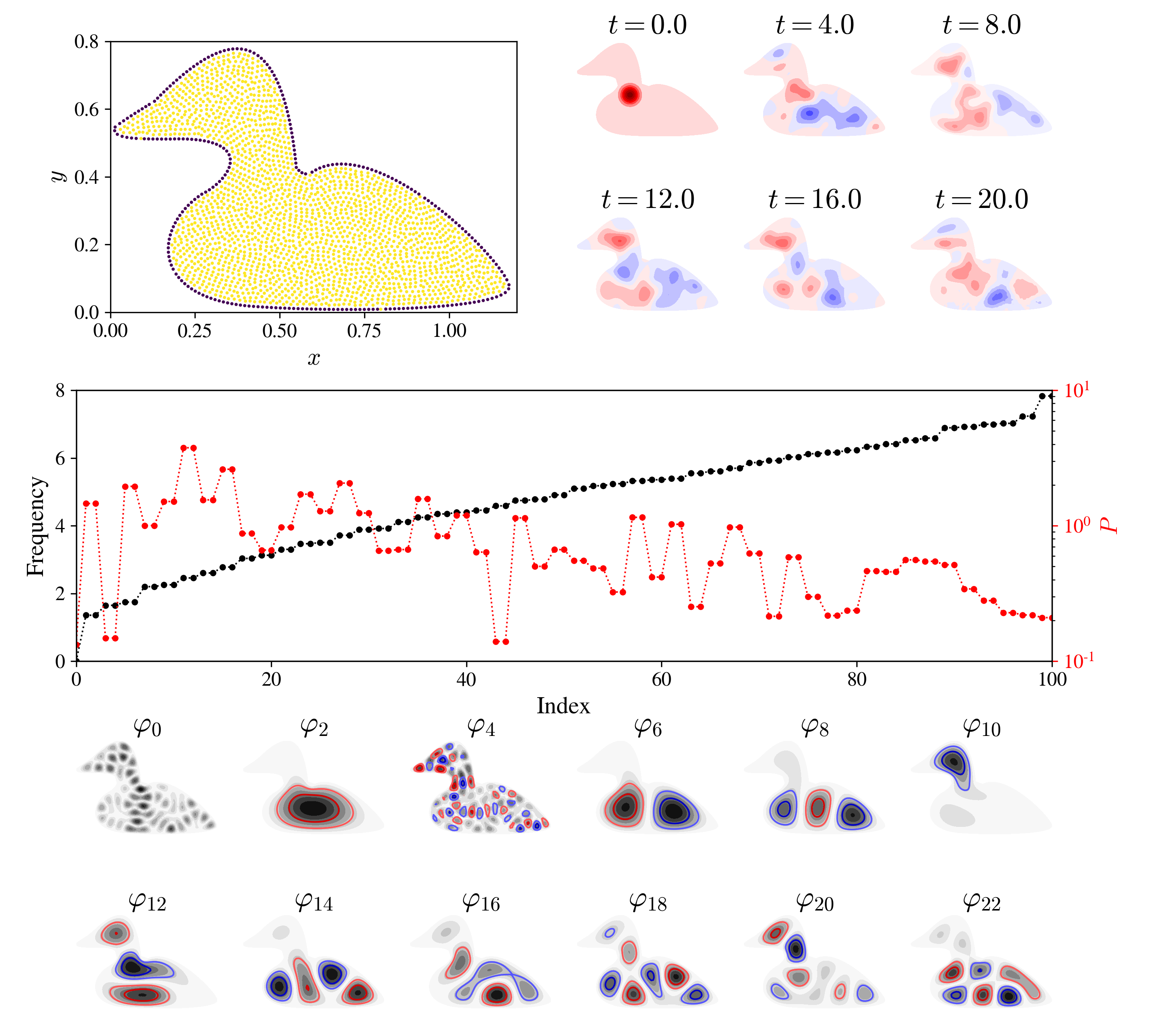}
	\caption{The DMD decomposition applied to a duck-shaped oscillating 
	membrane. \textit{top left:} Positioning of computational nodes with yellow 
	for interior nodes where we solve the wave equation and purple for boundary 
	nodes with enforced Dirichlet condition. \textit{top right:} A subset of 
	membrane state snapshots used for DMD. \textit{centre:} The DMD 
	spectrum showing the frequencies and powers of modes ordered by the former. 
	\textit{bottom:} Eigenvectors of DMD modes ordered by frequency.}
	\label{fig:duckMembrane}
\end{figure}

Now that we have shown that DMD can be used to identify inherent oscillatory 
modes we apply it to an irregularly shaped membrane. The duck-shaped membrane 
is again disturbed with a Gaussian initial state, resulting in the oscillation 
shown in Figure~\ref{fig:duckMembrane}. This combination of meshless 
computational nodes, the time evolution of the system, the DMD spectrum, and 
the DMD eigenvectors ordered by frequency provides a condensed overview of the 
dynamics and shows how DMD could be used in engineering analysis. The DMD 
spectrum is a convenient and commonly used display of DMD mode frequency and 
power in the same graph, which can, when sorted by frequency, be used similarly 
to the Fourier transform.

\section{Non-Newtonian fluid}
\label{ch:DVD}

The last example is a practical use case for DMD applied to hydrodynamics. We 
solve, with details available in~\cite{rot2022refined}, the system of partial 
differential equations
\begin{align}
	\div \vec{v} &= 0, \label{eq:physics1}\\
	\rho (\pdv{\vec{v}}{t} + \vec{v} \cdot \grad{\vec{v}}) &= -\grad p + 
	\div(\eta \grad
	\vec{v}) -\vec{g} \rho \beta T_\Delta, \label{eq:physics2}\\
	\rho c_p (\pdv{T}{t} + \vec{v} \cdot \grad{T}) &= \div(\lambda \grad
	T),\label{eq:physics3}\\
	\eta &= \eta_0 \left(\frac{1}{2}\norm{\grad{\vec{v}} + 
		(\grad{\vec{v}})^T}
	\right) ^{\frac{n-1}{2}}.\label{eq:physics4}
\end{align}
with $\vec{v}$, $T$, $p$, $\rho$, $\vec{g}$, $\beta$, $T_\Delta$, $c_p$, 
$\eta_0$, $n$ representing the flow velocity field, temperature field, 
pressure field, density, gravity, thermal expansion coefficient, 
temperature offset, heat capacity, viscosity constant and non-Newtonian 
index respectively.

The system is used to describe the natural convection in an incompressible 
non-Newtonian fluid. Non-Newtonian fluids have a variable viscosity meaning 
that the relationship between the shear-strain and the shear-stress is 
non-linear. We use a simple power-law model to describe the relationship with 
non-Newtonian power index $n$ as the main parameter of behaviour. Smaller 
values of $n$ indicate stronger shear-thinning non-Newtonian effects while 
$n=1$ equals to a normal, Newtonian, fluid.

We will be using the dimensionless Rayleigh number (Ra) that has a similar 
meaning for natural convection as Reynolds number, meaning stronger dynamics 
for higher values, as a case parameter, with details unimportant for the 
demonstration.

\begin{figure}
	\includegraphics[width=\linewidth]{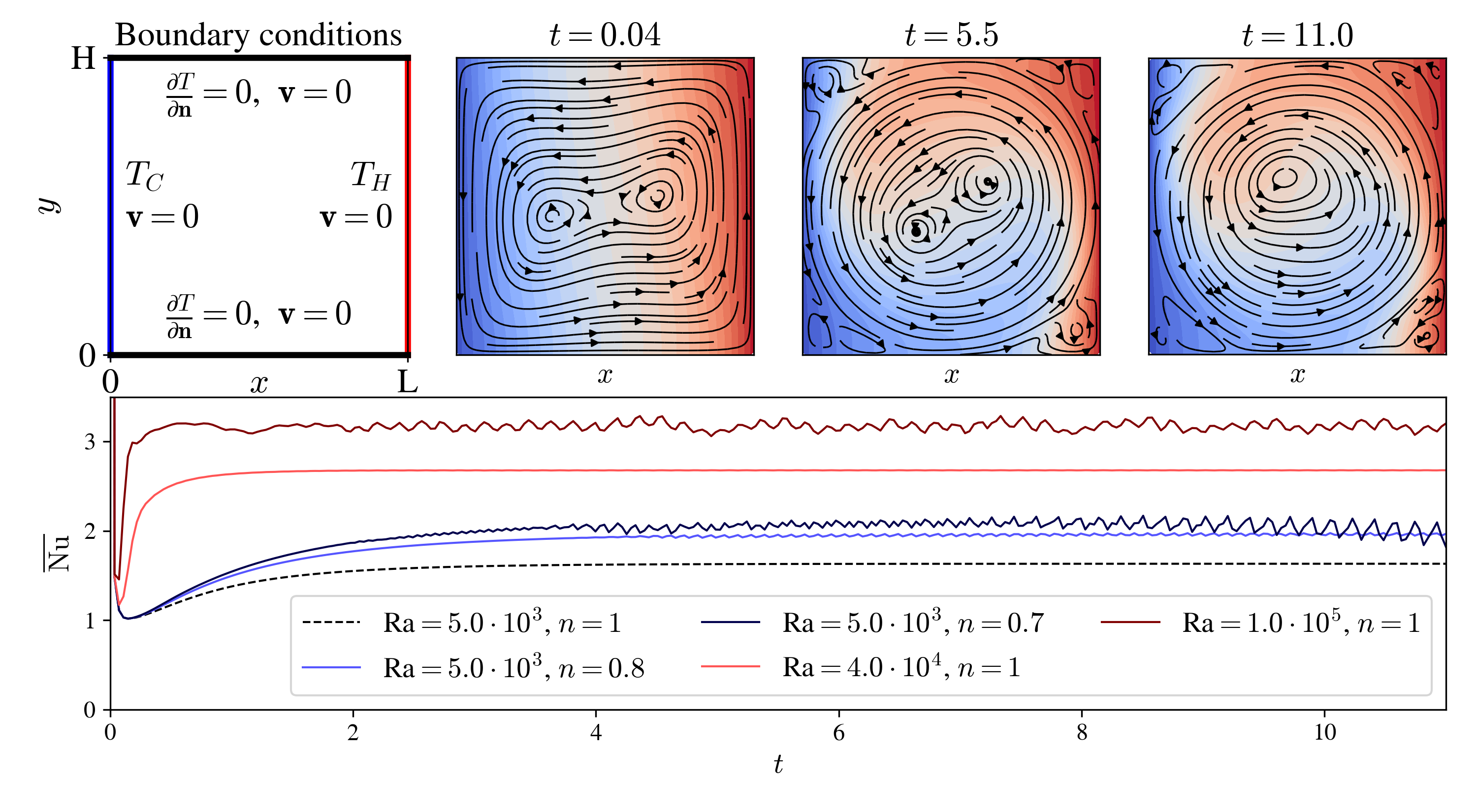}
	\caption{The setup and behaviour for the differentially heated 
		cavity non-Newtonian natural convection case. \textit{upper left: } The 
		case setup with boundary conditions. \textit{upper right:} A selection 
		of snapshots sampled from the system with Ra$=10^6$ and $n=1$. The 
		red-blue background colour shows the local fluid temperature with added 
		streamlines to illustrate the flow direction. \textit{bottom:} Time 
		evolution of Nusselt number for different parameters.}
	\label{fig:caseDVD}
\end{figure}

The model is solved in a square cavity with differentially heated walls 
schematically shown in the upper left graph of Figure~\ref{fig:caseDVD}. The 
left and right walls are kept at a temperature differential causing the fluid 
to form a vortex as it heats and rises at the right boundary and cools and 
descends on the left. This circulation is stationary for low Ra with a 
transition into an oscillatory and later chaotic regime as Ra increases. The 
upper right graphs in Figure~\ref{fig:caseDVD} show a selection of system 
temperature and velocity field snapshots for the oscillatory behaviour at 
Ra$=10^6$. A set of parameters where the transition into the oscillatory 
behaviour occurs has been identified by Kosec et al. \cite{Kosec2013} and we 
use DMD to analyse the differences in dynamics as we push past that point with 
increasing Ra and decreasing $n$.

As mentioned in the introduction we can use the Nusselt number, shown in the 
bottom graph of Figure~\ref{fig:caseDVD}, as an observable into the systems 
behaviour.  The system is stationary for Ra$=5 \cdot 10^3$ and $n=1$ but 
oscillations occur both when we increase Ra and when we decrease $n$.

\begin{figure*}
	\subfloat[Different Ra, Newtonian $n=1$.
		\label{fig:DVD_DMD_Ra}]
		{{\includegraphics[width=0.5\linewidth]{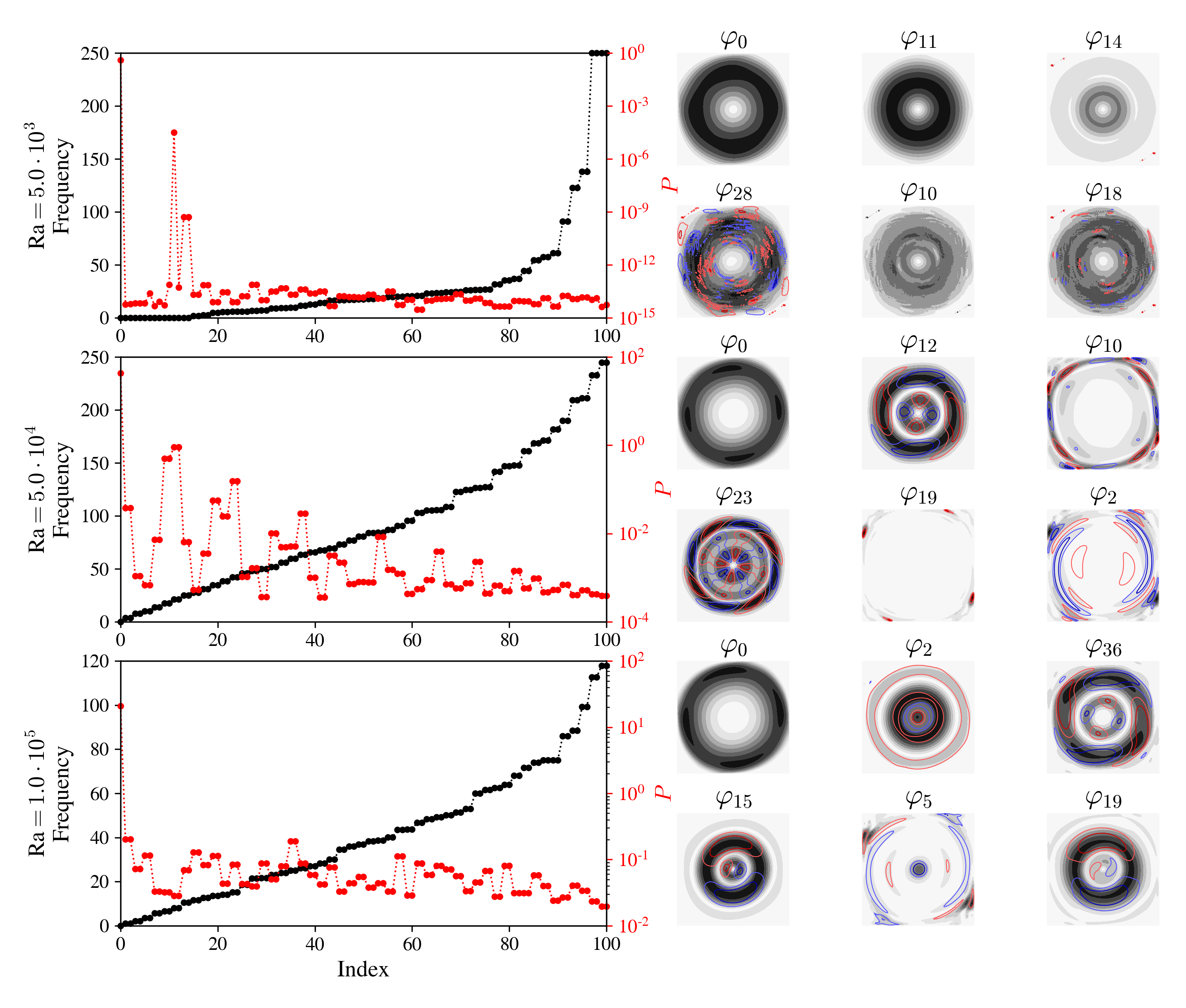} }}%
	\subfloat[Different non-Newtonian $n$, Ra$=5 \cdot 10^3$.
		\label{fig:DVD_DMD_n}]
		{{\includegraphics[width=0.5\linewidth]{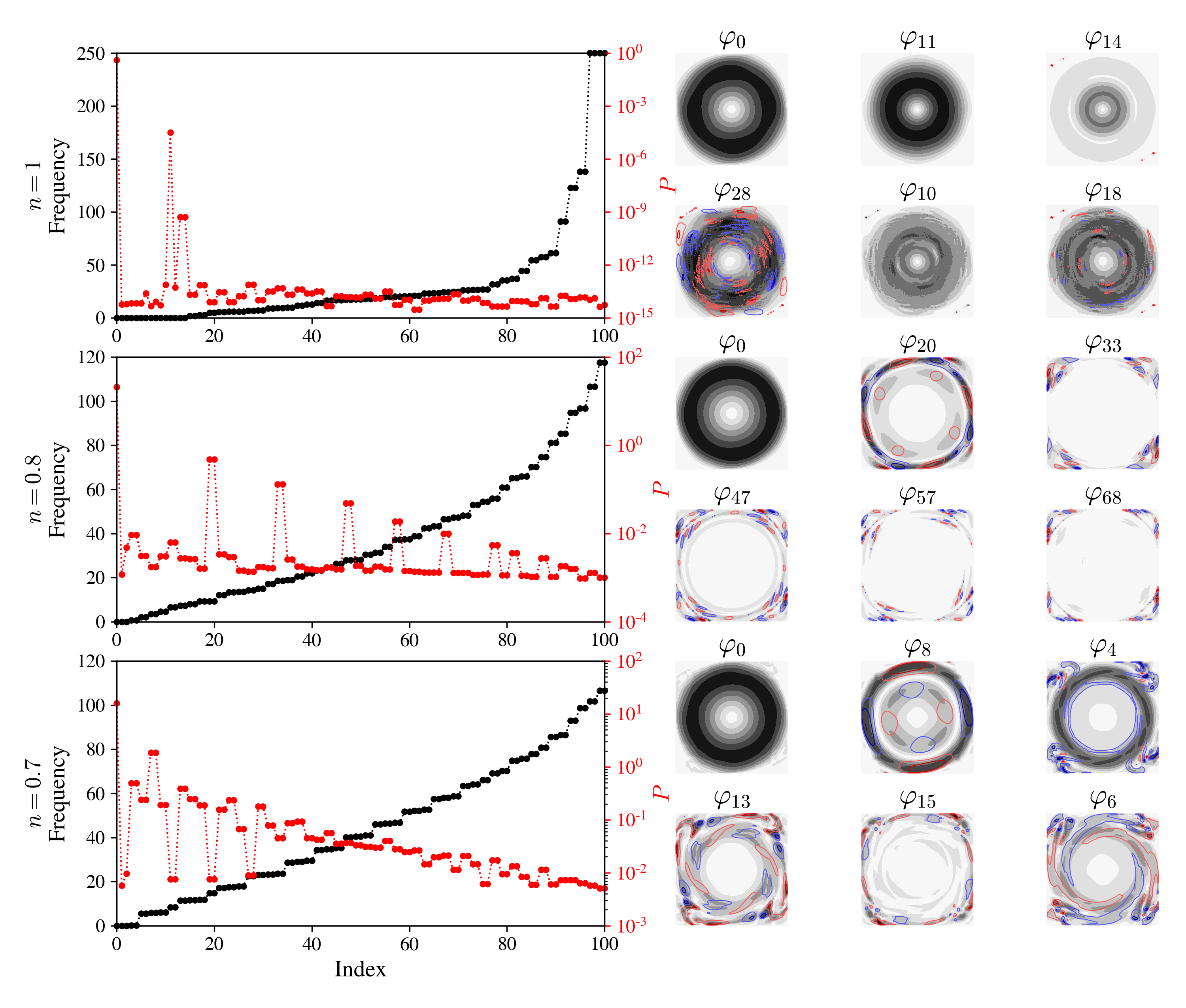} }}%
	\caption{The DMD spectra and eigenvectors for the strongest modes.}
\end{figure*}

We use DMD to decompose the 5 cases utilising 500 snapshots of the system 
between $t = 10$ and $t=11$ with stacking $m=5$. The results are shown in 
Figure~\ref{fig:DVD_DMD_Ra} for increasing Ra and Figure~\ref{fig:DVD_DMD_n} 
for decreasing $n$. The mode power defined in Eq.~\eqref{eq:modePower} is used 
to identify the strongest modes and display the corresponding eigenvectors. The 
first observation that is common to both modes of spurring the dynamics 
pertains to the DMD spectrum. The initial stationary case has the vast majority 
of its power in the first constant mode with others most likely only containing 
noise. As the intensity of the dynamics increases we get more distinct modes 
but the power of the background modes also increases as expected with ever 
wilder dynamics heading towards chaos, where all modes would be present and 
indistinguishable by power.

We can be reasonably certain that the oscillatory dynamics are in fact 
different for the Ra and $n$ spurred dynamics by using the eigenvectors of DMD 
modes sorted by power to identify the important parts. The strongest areas are 
clearly different with eigenvectors for the Ra induced dynamics in 
Figure~\ref{fig:DVD_DMD_Ra} present mainly as perturbations to the main central 
vortex while the modes for the $n$ induced dynamics in 
Figure~\ref{fig:DVD_DMD_n} point to activity in corners where counter vortices 
occur. This is consistent with the typical shear-thinning non-Newtonian 
behaviour where the larger velocity gradients are penalised less providing 
better conditions for smaller vortices.

\section{Conclusion}
We have presented the algorithm for dynamic mode decomposition and applied it 
to various test and example cases. The examples progressed from the most basic, 
where we were able to verify the results agains known closed-form values, 
towards a concrete hydrodynamics example where DMD provided valuable insight 
into the system that would be difficult to obtain otherwise. The DMD has proved 
to be a very practical tool for dynamical system analysis and it is recommended 
that anyone dealing with hydrodynamics or similarly complex systems should be 
at least familiar with its fundamentals. The algorithm is only a decade old and 
has an active community working on extensions for the concept of Koopman 
operator decomposition and applications for the existing versions of DMD.

\balance
\bibliographystyle{IEEEtran}
\bibliography{IEEEabrv, conference_dmd}

\begin{thebibliography}{10}
\providecommand{\url}[1]{#1}
\csname url@samestyle\endcsname
\providecommand{\newblock}{\relax}
\providecommand{\bibinfo}[2]{#2}
\providecommand{\BIBentrySTDinterwordspacing}{\spaceskip=0pt\relax}
\providecommand{\BIBentryALTinterwordstretchfactor}{4}
\providecommand{\BIBentryALTinterwordspacing}{\spaceskip=\fontdimen2\font plus
\BIBentryALTinterwordstretchfactor\fontdimen3\font minus
  \fontdimen4\font\relax}
\providecommand{\BIBforeignlanguage}[2]{{%
\expandafter\ifx\csname l@#1\endcsname\relax
\typeout{** WARNING: IEEEtran.bst: No hyphenation pattern has been}%
\typeout{** loaded for the language `#1'. Using the pattern for}%
\typeout{** the default language instead.}%
\else
\language=\csname l@#1\endcsname
\fi
#2}}
\providecommand{\BIBdecl}{\relax}
\BIBdecl

\bibitem{Ricciardi2019}
\BIBentryALTinterwordspacing
T.~R. Ricciardi, W.~Wolf, R.~L. Speth, and P.~Bent, ``Analysis of noise sources
  in realistic landing gear configurations through high fidelity simulations,''
  in \emph{AIAA Scitech 2019 Forum}, 2019. [Online]. Available:
  \url{https://arc.aiaa.org/doi/abs/10.2514/6.2019-0003}
\BIBentrySTDinterwordspacing

\bibitem{Lumley1967}
J.~L. Lumley, ``The structure of inhomogeneous turbulent flows,'' in
  \emph{Proceedings of the International Colloquium on the Fine Scale Structure
  of the Atmosphere and its Influence on Radio Wave Propagation}, A.~M. Yaglam
  and V.~I.Tatarsky, Eds., Moscow, 1967, pp. 166--178.

\bibitem{Rowley2005}
C.~W. Rowley, ``Model reduction for fluids, using balanced proper orthogonal
  decomposition,'' \emph{International Journal of Bifurcation and Chaos},
  vol.~15, no.~03, pp. 997--1013, 2005.

\bibitem{Taira2017}
\BIBentryALTinterwordspacing
K.~Taira, S.~L. Brunton, S.~T.~M. Dawson, C.~W. Rowley, T.~Colonius, B.~J.
  McKeon, O.~T. Schmidt, S.~Gordeyev, V.~Theofilis, and L.~S. Ukeiley, ``Modal
  analysis of fluid flows: An overview,'' \emph{AIAA Journal}, vol.~55, no.~12,
  pp. 4013--4041, 2017. [Online]. Available:
  \url{https://doi.org/10.2514/1.J056060}
\BIBentrySTDinterwordspacing

\bibitem{Holmes2012}
P.~Holmes, J.~L. Lumley, G.~Berkooz, and C.~W. Rowley, \emph{Turbulence,
  Coherent Structures, Dynamical Systems and Symmetry}, 2nd~ed., ser. Cambridge
  Monographs on Mechanics.\hskip 1em plus 0.5em minus 0.4em\relax Cambridge
  University Press, 2012.

\bibitem{Schmid2010}
P.~J. Schmid, ``Dynamic mode decomposition of numerical and experimental
  data,'' \emph{Journal of Fluid Mechanics}, vol. 656, p. 5–28, 2010.

\bibitem{Proctor2015}
J.~Proctor and P.~Welkhoff, ``Discovering dynamic patterns from infectious
  disease data using dynamic mode decomposition,'' \emph{International health},
  vol.~7, pp. 139--45, 03 2015.

\bibitem{Grosek2014}
J.~Grosek and J.~N. Kutz, ``Dynamic mode decomposition for real-time
  background/foreground separation in video,'' \emph{arXiv e-prints}, p.
  arXiv:1404.7592, Apr. 2014.

\bibitem{Rowley2009}
C.~W. Rowley, I.~Mezić, S.~Bagheri, P.~Schlatter, and D.~S. Henningson,
  ``Spectral analysis of nonlinear flows,'' \emph{Journal of Fluid Mechanics},
  vol. 641, p. 115–127, 2009.

\bibitem{Tu2014}
J.~H. Tu, C.~W. Rowley, D.~M. Luchtenburg, S.~L. Brunton, and J.~N. Kutz, ``On
  dynamic mode decomposition: Theory and applications,'' \emph{Journal of
  Computational Dynamics}, vol.~1, no.~2, pp. 391--421, 2014.

\bibitem{Tu2013thesis}
J.~Tu, ``Dynamic mode decomposition: Theory and applications,'' Ph.D.
  dissertation, Princeton University, sep 2013.

\bibitem{Kutz2016}
\BIBentryALTinterwordspacing
J.~N. Kutz, S.~L. Brunton, B.~W. Brunton, and J.~L. Proctor, \emph{Dynamic Mode
  Decomposition: Data-Driven Modeling of Complex Systems}, ser. Other Titles in
  Applied Mathematics.\hskip 1em plus 0.5em minus 0.4em\relax Society for
  Industrial and Applied Mathematics, 2016. [Online]. Available:
  \url{https://books.google.si/books?id=wbedDAEACAAJ}
\BIBentrySTDinterwordspacing

\bibitem{tolstykh2003rbffd}
\BIBentryALTinterwordspacing
A.~I. Tolstykh and D.~A. Shirobokov, ``On using radial basis functions in a
  ``finite difference mode'' with applications to elasticity problems,''
  \emph{Computational Mechanics}, vol.~33, no.~1, pp. 68--79, Dec 2003.
  [Online]. Available: \url{https://doi.org/10.1007/s00466-003-0501-9}
\BIBentrySTDinterwordspacing

\bibitem{medusa}
\BIBentryALTinterwordspacing
J.~Slak and G.~Kosec, ``Medusa: A c++ library for solving pdes using strong
  form mesh-free methods,'' \emph{ACM Trans. Math. Softw.}, vol.~47, no.~3,
  Jun. 2021. [Online]. Available: \url{https://doi.org/10.1145/3450966}
\BIBentrySTDinterwordspacing

\bibitem{rot2022refined}
\BIBentryALTinterwordspacing
M.~Rot and G.~Kosec, ``Refined rbf-fd analysis of non-newtonian natural
  convection,'' 2022. [Online]. Available:
  \url{https://arxiv.org/abs/2202.08095}
\BIBentrySTDinterwordspacing

\bibitem{Kosec2013}
G.~Kosec and B.~Šarler, ``Solution of a low prandtl number natural convection
  benchmark by a local meshless method,'' \emph{International Journal of
  Numerical Methods for Heat \& Fluid Flow}, vol.~23, p.~22, 01 2013.

\end{thebibliography}

\end{document}